\documentclass[twoside,a4paper,12pt,leqno]{amsart}

\usepackage{verbatim,amsxtra,amssymb,amsmath,color,psfrag,mathrsfs,graphicx}

\definecolor{mygrey}{gray}{0.70}
\definecolor{mygreen}{rgb}{0,.75,0}
\definecolor{myred}{rgb}{1,0,0}
\definecolor{orange}{rgb}{1,.5,0}

\numberwithin{equation}{part}

\theoremstyle{definition}
\newtheorem{defn}[equation]{Definition}

\theoremstyle{remark}
\newtheorem{rem}[equation]{Remark}
\newtheorem{ex}[equation]{Example}

\newcommand{\exref}[1]{Example~\ref{#1}}
\newcommand{\figref}[1]{Figure~\ref{#1}}
\newcommand{\secref}[1]{Section~{\bf\ref{#1}}}

\newcommand\A{\mathfrak A}
\renewcommand\a{\alpha}
\newcommand\Aut{\operatorname{\mathsf{Aut}}}
\newcommand\aut{\operatorname{\mathsf{Aut}}}

\renewcommand\b{\beta}
\newcommand\BB{{\mathcal B}}
\newcommand{\bded}{{\mathfrak{BA}}}
\newcommand{\bdd}{{\mathfrak B}}

\newcommand{\C}{\mathbb{C}}

\renewcommand{\d}{\delta}

\newcommand{\e}{\varepsilon}

\newcommand\F{\mathcal F}
\newcommand\f{\varphi}
\newcommand\fnst{{\mathfrak A}}

\newcommand\G{\Gamma}
\newcommand\GG{{\mathfrak G}}

\renewcommand\k{\kappa}

\renewcommand\l{\lambda}

\newcommand{\mapstoto}{\mathop{\,\mapstochar\relbar\joinrel\relbar\joinrel\longrightarrow\,}}
\newcommand\mthr{{\mathfrak M}}

\newcommand\ov{\overline}

\def\pair<#1>{{\langle\!\langle}#1{\rangle\!\rangle}}
\renewcommand\Pr{\mathcal P}
\newcommand\pt{\partial}

\renewcommand\S{{\mathsf S}}
\newcommand\s{\sigma}
\newcommand\sd{\rightthreetimes}

\newcommand\supp{\operatorname{\mathsf{supp}}}
\newcommand\sym{\operatorname{\mathsf{Sym}}}

\newcommand\toto{\mathop{\;\longrightarrow\;}}

\newcommand{\vn}{\varnothing}

\newcommand\wt{\widetilde}

\newcommand{\Z}{\mathbb Z}

\begin{document}

\title{Self-similarity and random walks}

\author{Vadim A. Kaimanovich}

\address{Jacobs University Bremen, 28759 Bremen, Germany}

\email{vadim.kaimanovich@gmail.com}

\begin{abstract}
This is an introductory level survey of some topics from a new branch of
fractal analysis --- the theory of self-similar groups. We discuss recent
works on random walks on self-similar groups and their applications to the
problem of amenability for these groups.
\end{abstract}

\maketitle

\thispagestyle{empty}

\section*{Introduction}

The purpose of this paper is to give a brief survey of some ideas and methods
associated with new progress in understanding the so-called \emph{self-similar
groups}. This class of groups consists of automorphisms of homogeneous rooted
trees defined in a simple recursive way. A good account of the initial period
of the theory can be found in the survey
\cite{Bartholdi-Grigorchuk-Nekrashevych03} by Bartholdi, Grigorchuk and
Nekrashevych and in the recent monograph of Nekrashevych \cite{Nekrashevych05}
(their authors are among the most active contributors to this field).

Self-similar groups have a natural interpretation in terms of fractal
geometry; their \emph{limit sets} are very interesting fractal sets (see the
recent papers \cite{Nekrashevych-Teplyaev08,Rogers-Teplyaev09} for a study of
Laplacians on limit sets). These groups also arise as \emph{iterated monodromy
groups} of rational endomorphisms of the Riemann sphere, which, for instance,
led to a recent solution of an old problem from rational dynamics
\cite{Bartholdi-Nekrashevych06}.

Self-similar groups are often quite unusual from the point of view of the
traditional group theory. This has both advantages and disadvantages. On one
hand, by using self-similar groups it is easy to construct examples which may
otherwise be much less accessible (for instance, the famous \emph{Grigorchuk
group of intermediate growth} \cite{Grigorchuk80,Grigorchuk85,delaHarpe00} has
a very simple self-similar presentation). On the other hand, even the simplest
group theoretical questions for self-similar groups may be quite hard. Already
finding a self-similar realization of a free group is very far from being
obvious \cite{Brunner-Sidki98,Glasner-Mozes05,Vorobets-Vorobets07}.

Another question of this kind is whether a given self-similar group is
\emph{amenable} (ame\-na\-bi\-lity introduced by von Neumann
\cite{vonNeumann29} is, in a sense, the most natural generalization of
finiteness, and it plays a fundamental role in group theory). There are
numerous characterizations of amenability in various terms. In particular, it
is known to be equivalent to existence of a \emph{random walk} on the group
with \emph{trivial behaviour at infinity} ($\equiv$ \emph{trivial Poisson
boundary}) \cite{Furstenberg73, Rosenblatt81,Kaimanovich-Vershik83}. It turns
out that self-similar groups may have random walks which are also self-similar
in a certain sense, and it is this self-similarity that can be used in order
to prove triviality of the Poisson boundary, and therefore establish
amenability of the underlying group.

This idea was first used by Bartholdi and Vir\`ag \cite{Bartholdi-Virag05} for
proving amenability of the so-called \emph{Basilica group} $\BB$. This group
first studied by Grigorchuk and \.Zuk \cite{Grigorchuk-Zuk02a} has a very
simple matrix presentation and also arises as the iterated monodromy group of
the map $z\mapsto z^2-1$ (known as the Basilica map, whence the name). In
particular, Grigorchuk and \.Zuk proved that $\BB$ is not
\emph{subexponentially elementary}, which made especially interesting the
question about its amenability. The approach of Bartholdi and Vir\`ag was
further developed by Kaimanovich \cite{Kaimanovich05} who used the
\emph{entropy theory} of random walks (which provides a simple criterion of
triviality of the Poisson boundary) in combination with a contraction property
for the asymptotic entropy of random walks on self-similar groups (the
\emph{``M\"unchhausen trick''}).

An important ingredient of this technique is a link (established in
\cite{Kaimanovich05}) between self-similarity and the so-called \emph{random
walks with internal degrees of freedom (RWIDF)} \cite{Kramli-Szasz83} also
known under the names of \emph{matrix-valued random walks}
\cite{Connes-Woods89} or of \emph{covering Markov chains}
\cite{Kaimanovich95}. These are group invariant Markov chains which take place
on the product of a group by a certain parameter set. Any random walk on a
self-similar group naturally gives rise to a random walk with internal degrees
of freedom parameterized by the alphabet of the action tree of the group. In
turn, this RWIDF, when restricted to the copy of the group corresponding to a
fixed value of the freedom parameter, produces a new random walk on the
self-similar group (as it is pointed out in \cite{Grigorchuk-Nekrashevych07},
this transformation corresponds to the classical operation of taking the
\emph{Schur complement} of a matrix). It is the interplay between the original
and the new random walks, which allows one to apply the M\"unchhausen trick.

This technique was recently applied by Bartholdi, Kaimanovich and Nekrashevych
\cite{Bartholdi-Kaimanovich-Nekrashevych08} to prove amenability of all
\emph{self-similar groups generated by bounded automata} (this class, in
particular, contains the Basilica group).

\medskip

In this survey we attempt to give a historical and conceptual overview of
these developments without going into technical details, so that it should
hopefully be suitable for a first acquaintance with the subject. Structurally,
our presentation is split into three parts. In \secref{sec:1} we discuss the
notion of self-similarity, introduce self-similar groups in general and the
subclass of self-similar groups generated by bounded automata. Further in
\secref{sec:2} we briefly discuss the notion of amenability of a group.
Finally, in \secref{sec:3} we analyze random walks on self-similar groups, and
show how they can be used for establishing amenability of self-similar groups.

The presentation is based on a talk given at the ``Fractal Geometry and
Stochastics IV" conference (and on several other occasions as well). I would
like to thank the organizers of this meeting for a very interesting and
inspiring week in Greifswald.

\section{Self-similar groups} \label{sec:1}

\subsection{And so \emph{ad infinitum} \dots}

The modern idea of self-similarity is best described by the following quote
from \emph{On Poetry: a Rhapsody} by Jonathan Swift (1733)\footnotemark :

\footnotetext{It was \emph{naturally extended} by Augustus de Morgan in his
\emph{Budget of Paradoxes} (1872):

\begin{quote}
Great fleas have little fleas upon their backs to bite 'em, \\
And little fleas have lesser fleas, and so \emph{ad infinitum.} \\
And the great fleas themselves, in turn have, greater fleas to go on; \\
While these again have greater still, and greater still, and so on.
\end{quote}

This is a description of what is nowadays called the \emph{natural extension}
of a non-invertible dynamical system. See \cite{Kaimanovich03a} for its
applications in the context of fractal sets.}

\begin{quote}
So, naturalists observe, a flea \\
Has smaller fleas that on him prey; \\
And these have smaller still to bite 'em, \\
And so proceed \emph{ad infinitum}.
\end{quote}

On a more formal level, the simplest self-similarity assumptions are:

\begin{itemize}
\item a part is similar to the whole;
\item the whole is a union of such parts;
\item these parts are pairwise disjoint.
\end{itemize}

\begin{figure}[h]
\begin{center}
     \psfrag{a}[cl][cl]{$a$}
     \psfrag{b}[cl][cl]{$b$}
     \psfrag{aa}[cl][cl]{$aa$}
     \psfrag{ab}[cl][cl]{$ab$}
     \psfrag{ba}[cl][cl]{$ba$}
     \psfrag{bb}[cl][cl]{$bb$}
     \psfrag{o}[cl][cl]{$\varnothing$}
     \psfrag{ta}[cl][cl]{$T_a$}
          \includegraphics[scale=.6]{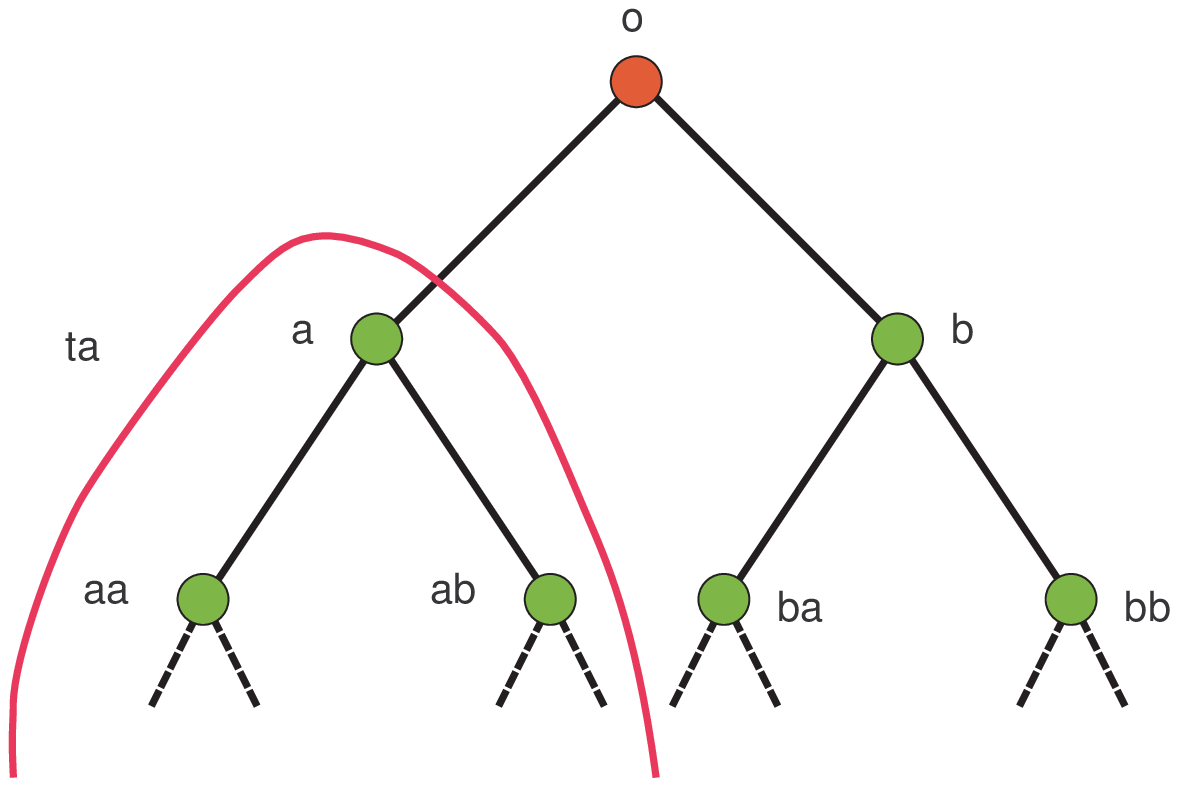}
          \end{center}
          \caption{}
          \label{fig:tree}
\end{figure}

These assumptions naturally lead to the most fundamental self-similar
structure, that of a \emph{rooted homogeneous tree}. Such trees are, for
instance, skeletons of iterated function systems satisfying the open set
condition (e.g., see \cite{Falconer03}), the simplest one of which is the
classical Cantor set. More precisely, let $X\cong\{1,2,\dots,d\}$ be a finite
set called the \emph{alphabet}. Denote by $X^*$ the set of all finite words in
the alphabet $X$ (including the empty word $\vn$). In other terms, $X^*$ is
the \emph{free monoid} generated by the set $X$ (the composition being the
concatenation $(w,w')\mapsto ww'$). The associated rooted homogeneous tree
$T=T(X)$ is the (right) Cayley graph of the free monoid $X^*$ (so that one
connects $w$ to $wx$ by an edge for all $w\in X^*,x\in X$). The tree
$T(X)\cong X^*$ is split into \emph{levels} $T_n\cong X^n$ (the set of words
of length $n$). The level $T_0$ consists only of the empty word $\vn$, which
is the root of $T$. Each vertex $w\in T\cong X^*$ is the root of the subtree
$T_w$ which consists of all the words beginning with $w$. The map $w'\mapsto
ww'$ provides then a canonical identification of the trees $T$ and $T_w$, see
\figref{fig:tree}, where $X=\{a,b\}$.

\subsection{Generalized permutation matrices} \label{sec:perm}

\begin{figure}[h]
\begin{center}
     \psfrag{t}[cl][cl]{$T_x$}
     \psfrag{x}[cl][cl]{$x$}
     \psfrag{x1}[cl][cl]{$x'$}
     \psfrag{gx}[cl][cl]{$g_x$}
     \psfrag{tx}[cl][cl]{$T_{x'}$}
\includegraphics[scale=.6]{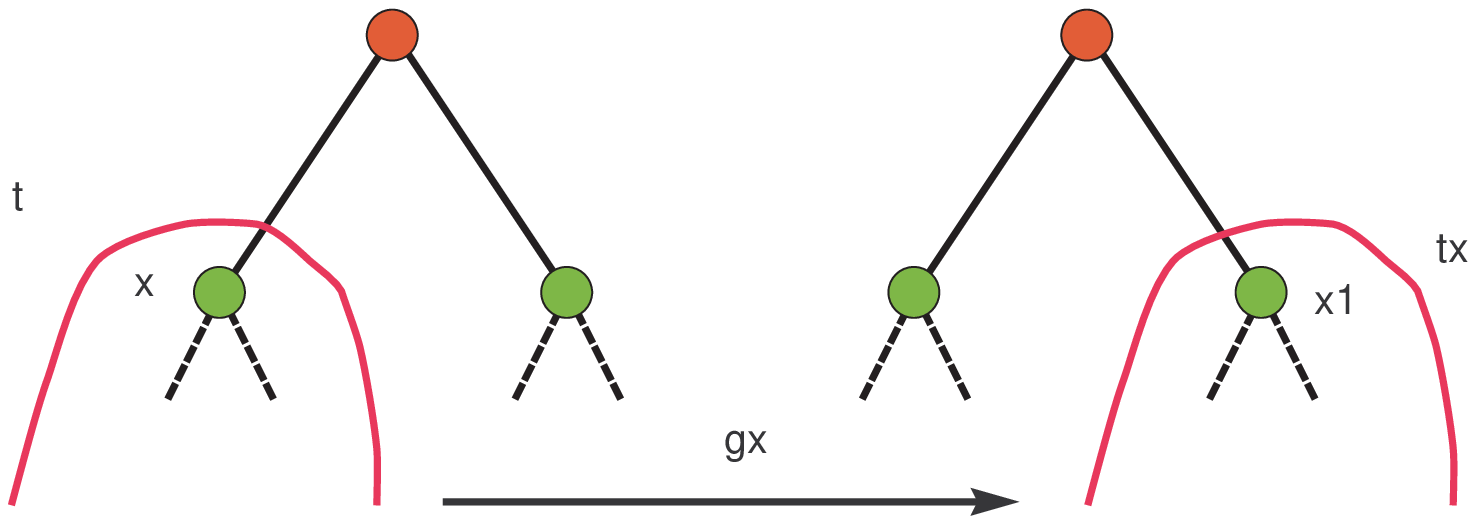}
\end{center}
\caption{}
\label{fig:two}
\end{figure}

The group $\GG=\Aut(T)$ of automorphisms of the tree $T$ obviously preserves
each level of $T$. In particular, it acts by permutations on the first level
$X\cong T_1$, i.e., there is a homomorphism $g\mapsto\s=\s^g$  from $\GG$ to
$\sym(X)$ (the permutation group on $X$). Another piece of data associated
with any automorphism $g\in\GG$ is a collection of automorphisms
$\{g_x\}_{x\in X}$ indexed by the alphabet $X$. Indeed, if $y=\s^g(x)$, then
$g$ establishes a one-to-one correspondence between the corresponding subtrees
$T_x$ and $T_y$ rooted at the points $x$ and~$y$, respectively. Since both
subtrees $T_x$ and $T_y$ are canonically isomorphic to the full tree $T$, the
map $g:T_x\to T_y$ is conjugate to an automorphism of $T$ denoted by $g_x$ (in
terms of Swift's description above, $g_x$ describes what happens on the back
of a first order flea when it moves from position $x$ to position
$y=\s^g(g)$), see \figref{fig:two}.

Conversely, a permutation $\s\in\sym(X)$ and a collection $\{g_x\}_{x\in X}$
of elements of $\GG$ uniquely determine the associated automorphism $g\in\GG$.
In algebraic terms it means that $g\mapsto \bigl(\s^g;\{g_x\}_{x\in X}\bigr)$
is an isomorphism of the group $\GG$ and the \emph{semi-direct product}
$\sym(X)\sd\GG^X$ (in yet another terminology, $\GG$ is isomorphic to the
\emph{permutational wreath product} $\GG\wr\sym(X)$). There is a very
convenient way of visualizing this structure by means of \emph{generalized
permutation matrices.}

\begin{figure}[h]
\begin{center}
 \psfrag{a}[cl][cl]{$g_1$}
 \psfrag{b}[cl][cl]{$g_2$}
\includegraphics[scale=.6]{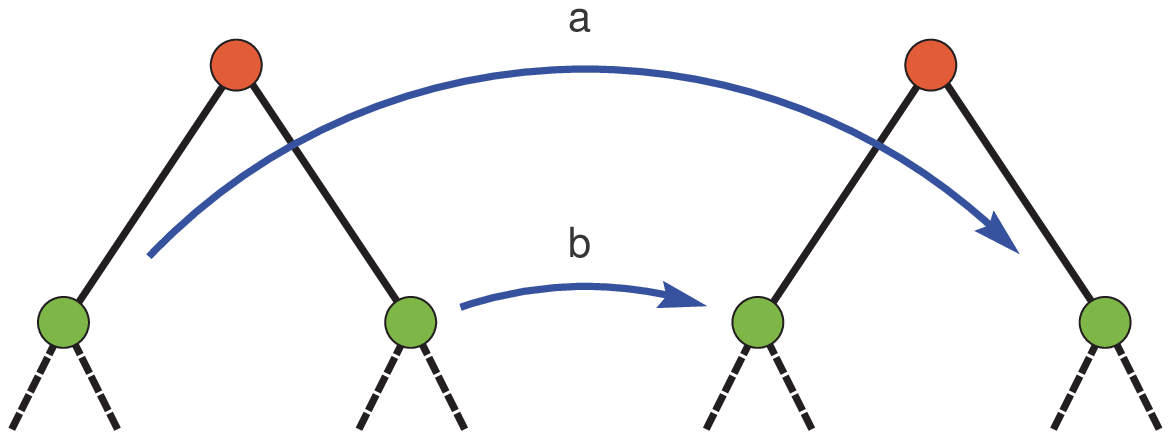}
\end{center}
\caption{} \label{fig:g12}
\end{figure}

Recall that the usual \emph{permutation matrix} $M^\s$ associated with a
permutation $\s\in\sym(X)$ is a $|X|\times|X|$ matrix with entries
$$
M^\s_{xy}=\left\{\begin{array}{ll}1 \;, &\text{if $y=\sigma(x)$,}\\
0 \;, & \text{otherwise}\end{array}\right. \;,
$$
and that the map $\s\mapsto M^\s$ is a group isomorphism. In the same way we
shall present the data $\bigl(\s^g;\{g_x\}_{x\in X}\bigr)$ by the
\emph{generalized permutation matrix} $M^g$ with entries
$$
M^g_{xy}=\left\{\begin{array}{ll}g_x \;, &\text{if $y=\sigma^g(x)$,}\\
0 \;, & \text{otherwise.}\end{array}\right.
$$
For instance, the automorphism $g$ described in \figref{fig:g12} is presented
by the matrix
$$
M^g = \begin{pmatrix} 0 & g_1 \\ g_2 & 0 \end{pmatrix} \;.
$$

More generally, given an arbitrary group $G$, we shall denote by
$$
\sym(X;G) = G\wr\sym(X)=\sym(X)\sd G^X
$$
the \emph{group of generalized permutation matrices} of order $|X|$ with
non-zero entries from the group $G$. The group operation here is the usual
matrix multiplication, the only difference with ordinary permutation matrices
being that the matrix elements are multiplied according to the group law of
$G$. Obviously, application of the \emph{augmentation map} (which consists in
replacing all group elements with~1) to a generalized permutation matrix
yields a usual permutation matrix, which corresponds to the natural
homomorphism of $\sym(X;G)$ onto $\sym(X)$. We can now sum up the above
discussion by saying that \emph{there is a natural isomorphism of the group
$\GG=\aut(T)$ of automorphisms of the tree $T=T(X)$ and of the generalized
permutation group $\sym(X;\GG)$}. It is this isomorphism that embodies the
self-similarity properties of the group $\GG$.

\subsection{Self-similar groups and matrix presentations} \label{sec:matr}

\begin{defn}
A countable subgroup $G\subset\GG$ is \emph{self-similar} if the restriction
of the isomorphism $\GG\to\sym(X;\GG)$ to $G$ induces an embedding
$G\hookrightarrow \sym(X;G)$; in other words, if all entries of the matrices
$M^g,\,g\in G$ belong to $G$. Note that, rigorously speaking, self-similarity
is a property of the \emph{embedding} $G\subset\GG$ rather than of the group
$G$ only. The embedding $G\hookrightarrow \sym(X;G)$ need \emph{not} be
surjective (see the example below).
\end{defn}

\begin{ex}
The \emph{adding machine} (isomorphic to the group $\Z$) is generated by the
transformation $a:z\mapsto z+1$ on the ring $\Z_2$ of \emph{2-adic integers}
$\e_0 + \e_1\cdot 2 + \dots + \e_n\cdot 2^n + \dots$, where the digits $\e_i$
take values 0 or 1. Depending on the values of initial digits, it acts in the
following way:
$$
\begin{matrix}
  0 \!\!\!\!& \e_1 \!\!\!\!& \e_2 \!\!\!\!& \e_3  \dots \quad \mapsto \quad \!\!\!\!& 1 \!\!\!\!& \e_1 \!\!\!\!& \e_2 \!\!\!\!& \e_3  \dots \;, \\
  1 \!\!\!\!& 0    \!\!\!\!& \e_2 \!\!\!\!& \e_3  \dots \quad \mapsto \quad \!\!\!\!& 0 \!\!\!\!& 1    \!\!\!\!& \e_2 \!\!\!\!& \e_3  \dots \;, \\
  1 \!\!\!\!& 1    \!\!\!\!& 0    \!\!\!\!& \e_3  \dots \quad \mapsto \quad \!\!\!\!& 0 \!\!\!\!& 0    \!\!\!\!& 1    \!\!\!\!& \e_3  \dots \;, \\
 &&&  \phantom{\e_3\dots} \quad \hfil \hbox to 8mm{\dotfill} &&&
\end{matrix}
$$
We can think of sequences $(\e_0,\e_1,\dots)$ as of boundary points of the
binary rooted tree $T=T(X)$ of the alphabet $X=\{0,1\}$. The transformation
$a$ extends to an automorphism of $T$, and, as one can easily see from its
symbolic description above, the associated generalized permutation matrix is
$$
M^a = \begin{pmatrix} 0 & 1 \\ a & 0 \end{pmatrix} \;.
$$
Thus, the infinite cyclic group $\langle a \rangle\cong\Z$ generated by the
transformation $a$ is self-similar (as a subgroup of the full group of
automorphisms $\GG$).

Note that the automorphism $a$ is completely determined by the matrix $M^a$.
Indeed, the augmentation map applied to $M^a$ produces the usual permutation
matrix $\begin{pmatrix} 0 & 1 \\ 1 & 0 \end{pmatrix}$ which describes the
permutation by which $a$ acts on the first level of the tree $T$. Further, by
substituting $M^a$ for $a$ and the identity matrix for 1 in $M^a$ one obtains
the matrix
$$
\begin{pmatrix}
0 & 0 & 1 & 0 \\ 0 & 0 & 0 & 1 \\ 0 & 1 & 0 & 0 \\ a & 0 & 0 & 0
\end{pmatrix} \;,
$$
augmentation of which produces the order 4 permutation matrix
$$
\begin{pmatrix}
0 & 0 & 1 & 0 \\ 0 & 0 & 0 & 1 \\ 0 & 1 & 0 & 0 \\ 1 & 0 & 0 & 0
\end{pmatrix}
$$
describing the action of $a$ on the second level of $T$. By recursively
repeating this procedure one obtains the action of $a$ on all levels of $T$,
i.e., an automorphism of the whole tree $T$.
\end{ex}

The example above suggests the following way of defining self-similar groups
by their matrix presentations. Fix a finite (or countable, for infinitely
generated groups) set $K$ (the future set of generators of a self-similar
group), and assign to each $\k\in K$ a generalized permutation matrix $M^\k$
whose non-zero entries are words in the alphabet consisting of letters from
$K$ and their inverses (the entry associated with the empty word is 1). By
replacing the non-zero entries of matrices $M^\k$ with corresponding products
of the associated matrices and their inverses we obtain generalized
permutation matrices of order $|X|^2$, etc. The usual permutation matrices
obtained from them by the augmentation map determine then the action of
elements from $K$ on all levels of the tree $T$, i.e., the corresponding
automorphisms of $T$. See \cite{Bartholdi-Grigorchuk00} or
\cite{Nekrashevych05} for more on recursions of this kind.

A particular case of this construction arises in the situation when all
non-zero entries of matrices $M^\k$ are elements of the set $K$. In this case
the assignment $\k\mapsto M^\k$ amounts to a map $(\k,x)\mapsto (\l,y)$ of the
product $K\times X$ to itself, i.e., to an \emph{automaton}. Here $y=\s(x)$
for the permutation $\s=\s(M^\k)$ determined by the matrix $M^\k$, and $\l$ is
the matrix entry $M^\k_{xy}$. The self-similar group obtained in this way is
called an \emph{automaton group.}\footnote{Usually, when talking about
automata groups, one tacitly assumes that the corresponding automaton is
finite, i.e., the generating set $K$ is finite.}

\subsection{The Basilica group} \label{sec:bas}

The \emph{Basilica group} $\BB$ is determined by the matrix presentation
$$
 a \mapsto \begin{pmatrix} b & 0 \\ 0 & 1 \end{pmatrix} \;, \qquad
 b \mapsto \begin{pmatrix} 0 & a \\ 1 & 0 \end{pmatrix} \;.
$$
The aforementioned recursion for this group looks in the following way:
$$
\begin{aligned}
&a \mapsto \begin{pmatrix} b & 0 \\ 0 & 1 \end{pmatrix} \mapsto
\begin{pmatrix}
0 & a & 0 & 0 \\ 1 & 0 & 0 & 0 \\ 0 & 0 & 1 & 0 \\ 0 & 0 & 0 & 1
\end{pmatrix} \mapsto \dots \;, \\
&b \mapsto \begin{pmatrix} 0 & a \\ 1 & 0 \end{pmatrix} \mapsto
\begin{pmatrix}
0 & 0 & b & 0 \\ 0 & 0 & 0 & 1 \\ 1 & 0 & 0 & 0 \\ 0 & 1 & 0 & 0
\end{pmatrix} \mapsto \dots \;,
\end{aligned}
$$
and the associated automaton is
$$
\begin{cases}
 (a,1)&\mapsto\; (b,1) \\
 (a,2)&\mapsto\; (e,2) \\
 (b,1)&\mapsto\; (a,2) \\
 (b,2)&\mapsto\; (e,1) \\
 (e,1)&\mapsto\; (e,1) \\
 (e,2)&\mapsto\; (e,2) \;.
\end{cases}
$$
Here $\{1,2\}=X$ is the alphabet of the binary tree $T$, and $K=\{a,b,e\}$,
where $e$ is the group identity determined by the substitution $e\mapsto
\begin{pmatrix} e & 0\\ 0 & e \end{pmatrix}$.

The group $\BB$ was first studied by Grigorchuk and \.Zuk
\cite{Grigorchuk-Zuk02a} (see below for its algebraic propereties). The name
\emph{Basilica} comes from the fact that it also appears as the \emph{iterated
monodromy group} of the rational map $z\mapsto z^2-1$ on the Riemann sphere
$\ov\C$.

This latter notion was introduced by Nekrashevych who created a very fruitful
link between the theory of self-similar groups and \emph{rational dynamics}.
Namely, given a rational map $\f:\ov\C\to\ov\C$ of degree $d$, a generic point
$z\in\ov\C$ has precisely $d$ preimages, each of which also has $d$ preimages,
etc. Thus, attached to a generic point $z\in\ov\C$ is the rooted homogeneous
tree $T_z$ of its preimages. One can move the preimage tree along any
continuous curve consisting of generic points. However, if $z$ follows a
non-contractible loop, it may happen that, although the preimage tree $T_z$
returns to its original position, it undergoes a certain non-trivial
\emph{monodromy transformation}. Thus, there is a homomorphism of the
fundamental group of the connected component of $z$ in the set of generic
points to the group of automorphisms of the preimage tree $T_z$. The resulting
subgroup of $\aut(T_z)$ is called the \emph{iterated monodromy group} of the
map $\f$, see \cite{Nekrashevych05} for more details.

\begin{figure}[h]
\begin{center}
\includegraphics[scale=.4]{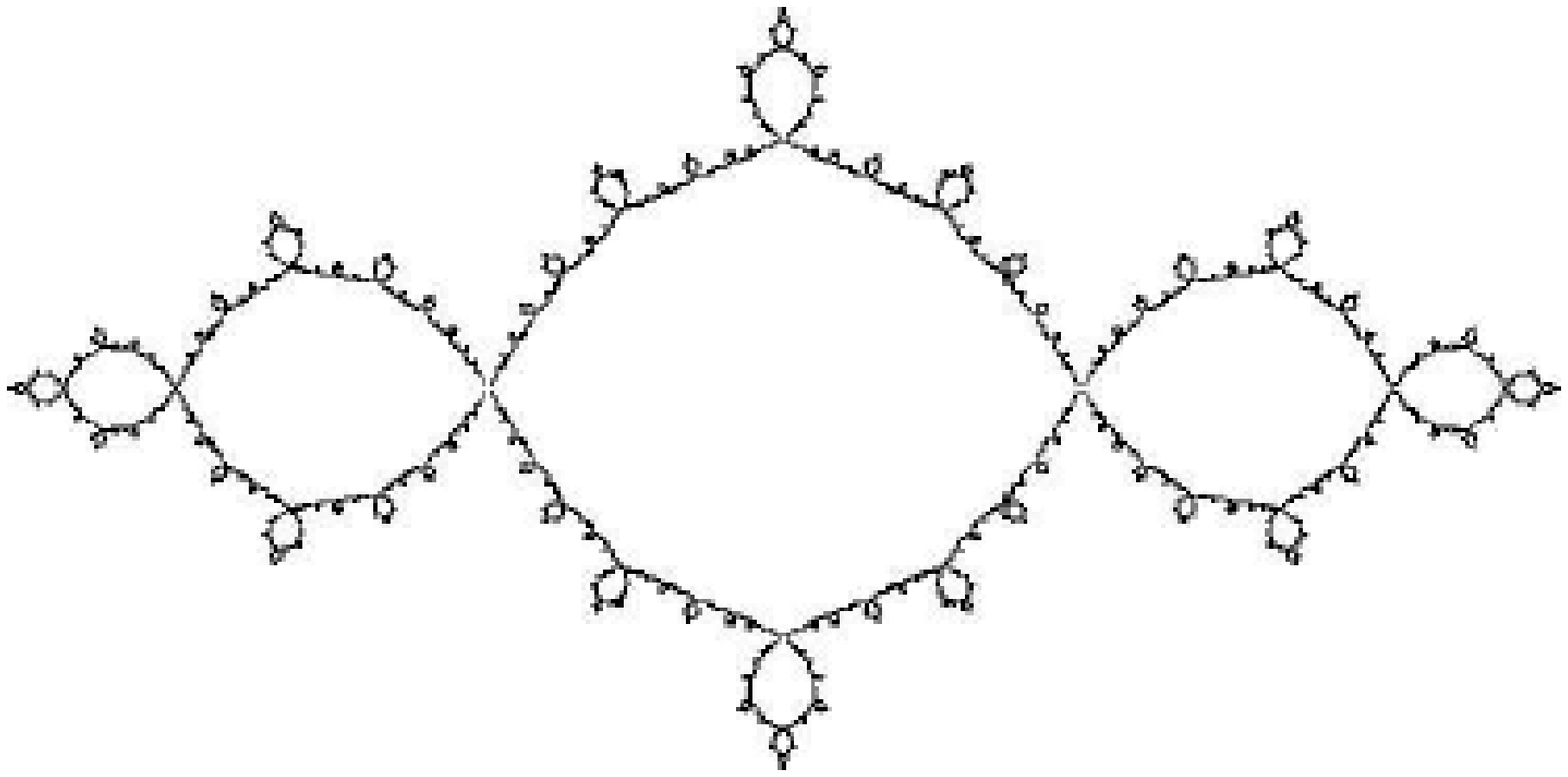}
\end{center}
\caption{} \label{fig:basilica}
\end{figure}

Now, in rational dynamics the map $z\mapsto z^2-1$ is called the
\emph{Basilica map} \cite{Bielefeld90}, because its \emph{Julia set} (a subset
of the Riemann sphere which, in a sense, consists of the limit points of this
map) looks similar to \emph{Basilica di San Marco} in Venice (together with
its reflection in the water), see \figref{fig:basilica}. This Julia set also
arises as the limit set of the Basilica group $\BB$.

Interesting algebraic and analytic properties of the Basilica group obtained
in \cite{Grigorchuk-Zuk02a, Grigorchuk-Zuk02b, Bartholdi-Grigorchuk02} made
especially relevant the question about its \emph{amenability} formulated in
\cite{Grigorchuk-Zuk02a}, also see \cite{Bartholdi-Grigorchuk-Nekrashevych03}.

\subsection{Bounded automatic automorphisms and Mother groups}
\label{sec:moth}

The Basilica group $\BB$ actually belongs to a certain natural subclass of the
class of self-similar groups which we shall now describe.

As we have seen in \secref{sec:perm}, given an automorphism $g\in\GG=\aut(T)$,
any symbol $x\in X$ determines an associated automorphism $g_x\in\GG$. In the
same way such an automorphism $g_w\in\GG$ (the \emph{state} of $g$ at the
point $w$) can be defined for an arbitrary word $w\in T\cong X^*$, by
restricting the automorphism $g$ to the subtree $T_w$ with the subsequent
identification of both $T_w$ and its image $g(T_w)=T_{g(w)}$ with $T$.
Equivalently, $g_w$ can be obtained by recursively applying the presentation
$g\mapsto M^g$, namely, $g_w$ is the non-zero entry of the row $R_w$ of the
$|X|^{|w|}\times|X|^{|w|}$ matrix obtained by recursive expansion $g\mapsto
M^g\mapsto \dots$.

If the \emph{set of states} of $g$
$$
\S(g)=\{g_w:w\in T\}\subset\GG
$$
is finite, then the automorphism $g$ is called \emph{automatic}. The set of
all automatic automorphisms of the tree $T$ forms a countable subgroup
$\fnst=\fnst(X)$ of $\GG=\GG(X)$.

An automorphism $g$ is called \emph{bounded} if the sets $\{w\in X^n:g_w\neq
e\}$ have uniformly bounded cardinalities over all $n$. The set of all bounded
automorphisms also forms a subgroup $\bdd=\bdd(X)$ of $\GG=\GG(X)$. We denote
by $\bded=\bded(X)=\bdd(X)\cap\fnst(X)$ the \emph{group of all bounded
automatic automorphisms} of the homogeneous rooted tree $T$.

\medskip

It is easy to see that the generators $a,b$ of the Basilica group $\BB$
described in \secref{sec:bas} are both automatic and bounded in the above
sense, so that $\BB\subset\bded$ (here and on several occasions below we omit
the alphabet $X$ from our notation). More generally, any group generated by a
(finite) automaton (see \secref{sec:matr}) is a subgroup of $\fnst$. An
automaton is called \emph{bounded} if its group is contained in $\bdd$
(therefore, in $\bded$). The class of \emph{groups generated by bounded
automata} was defined by Sidki in~\cite{Sidki00}. Obviously, all these groups
are subgroups of $\bded$. Most of the well-studied groups of finite automata
belong to this class (see \cite{Bartholdi-Kaimanovich-Nekrashevych08} for
examples).

Groups generated by bounded automata also appear naturally in connection with
fractal geometry. It was proved in~\cite{Bondarenko-Nekrashevych03} that every
such group is \emph{contracting}, and a contracting group is generated by
bounded automata if and only if the boundary of its \emph{tile} is finite.
This technical condition implies that the \emph{limit space}
(see~\cite{Nekrashevych05}) of such a group belongs to the well studied class
of \emph{nested fractals} (see~\cite{Kigami02}). This is the class of fractals
on which the Brownian motion is best understood. This shows an interesting
connection between the most well-studied class of self-similar groups and the
class of fractals with most well-understood analysis (see
\cite{Nekrashevych-Teplyaev08} for more details).

\medskip

It turns out that the class of groups generated by bounded automata contains a
countable family of groups which have certain \emph{universality properties}
with respect to this class.

Let $X$ be a finite set with a distinguished element $o\in X$, and put $\ov
X=X\setminus\{o\}$. Set $A=\sym(X)$ and $B=\sym(\ov X;A)$, and recursively
embed the groups $A$ and $B$ into $\GG(X)$ by the matrix presentations
$$
M^a = \phi_A(a) \;, \quad M^b=\begin{pmatrix} b & 0 \\ 0 & \phi_B(b)
\end{pmatrix} \;,
$$
where $\phi_A(a),\phi_B(b)$ are, respectively, the permutation and the
generalized permutation matrices corresponding to $a\in A,b\in B$. Then the
\emph{Mother group} $\mthr=\mthr(X)=\langle A, B\rangle$ is the subgroup of
$\GG=\GG(X)$ generated by the finite groups $A$ and $B$.

A direct verification shows that both groups $A,B$ are contained in $\bded$,
whence \emph{the group $\mthr(X)$ is a subgroup of $\bded(X)$}. On the other
hand, as it was proved in \cite{Bartholdi-Kaimanovich-Nekrashevych08},
\emph{any finitely generated subgroup of $\bded(X)$ can be embedded as a
subgroup into the generalized permutation group $\sym(X^N;\mthr(X^N))$ for
some integer $N$}.

\medskip

Thus, in view of the fact that amenability is preserved by \emph{elementary
operations} (see \secref{sec:elem} below), amenability of the groups
$\bded(X)$ for all finite sets $X$ (therefore, amenability of all groups
generated by bounded automata) would follow from amenability just of all the
Mother groups $\mthr(X)$.

It is worth noting that the groups generated by bounded automata form a
subclass of the class of \emph{contracting self-similar groups}
(see~\cite{Bondarenko-Nekrashevych03,Nekrashevych05}). It is still an open
question whether all contracting groups are amenable. However, Nekrashevych
\cite{Nekrashevych08} recently established a weaker property:
\emph{contracting groups contain no free groups with $\ge 2$ generators}.

\section{Amenability} \label{sec:2}

\subsection{From finite to infinite: von Neumann, Day and Reiter}

Finite groups can be characterized as those discrete groups which have a
\emph{finite invariant measure}. In other words, a discrete group $G$ is
finite if and only if the natural action of $G$ by translations on the space
$\ell^1_{+,1}(G)$ of positive normalized elements from $\ell^1(G)$ has a fixed
point. This trivial observation suggests two ways of ``extending'' the
finiteness property to infinite groups. One can look either for fixed points
in a bigger space, or for approximative invariance instead of exact one.

The first idea was implemented by John von Neumann \cite{vonNeumann29},
according to whose definition \emph{amenable groups are those which admit a
translation invariant mean}\footnotemark. A \emph{mean} on $G$ is a finitely
additive probability measure, in other words, an element of the space
$[\ell^\infty]^*_{+,1}$ of positive normalized functionals on
$\ell^\infty(G)$. Usual measures on $G$ are also means, but if $G$ is
infinite, then there are many more means than measures (which corresponds to
the fact that in the infinite case $\ell^1$ is significantly ``smaller'' than
its second dual space $[\ell^1]^{**}=[\ell^\infty]^*$).

\footnotetext{Actually, the original term used by von Neumann was the German
\emph{me{\ss}bare Gruppe}, which means ``measurable group'' in English. It was
later replaced in German with \emph{mittelbare} (cf. \emph{moyennable} in
French), literally meaning ``averageable''. In English, however, Mahlon M. Day
suggested to use (apparently, first as a pun) the word \emph{amenable}, which
appeared in print in this context for the first time in 1949 \cite{Day49}. It
is curious that Day himself, when he later described the history of this term
in \cite{Day83} on the occasion of the nomination of his paper \cite{Day57} as
a ``Citation Classic'' dated its appearance to 1955: \emph{In 1929, von
Neumann studied a new class of groups, those with invariant means on the
bounded functions. My thesis (1939) studied semigroups with invariant means;
thereafter, I worked in the field alternately with the geometry of Banach
spaces. I finished a large geometrical project in 1955 and turned back to
invariant means; in order to talk to my students I invented the term 'amenable
(pronounced as amean'able) semigroups'.}}

Means being highly non-constructive objects, the other way was explored
(surprisingly, much later than the original definition of von Neumann) by
Reiter \cite{Reiter65} who introduced (under the name $P_1$) what is nowadays
known as \emph{Reiter's condition} for a group $G$: \emph{there exists an
approximatively invariant sequence of probability measures on $G$,} in other
words, \emph{there exists a sequence of probability measures $\l_n$ on $G$
such that}
$$
\|\l_n - g\l_n\| \toto_{n\to\infty} 0 \qquad\forall\,g\in G \;,
$$
where $\|\cdot\|$ denotes the total variation norm. He proved that the above
condition is in fact equivalent to amenability as defined by von Neumann.

\begin{ex} \label{ex:cesaro}
The sequence of \emph{Cesaro averaging measures}
$$
\l_n = \frac1{n+1} \bigl (\d_0 + \d_1 + \dots + \d_n) \;,
$$
on the group of integers $\Z$ is approximatively invariant (here $\d_n$
denotes the unit mass at the point $n$). Thus, $\Z$ is amenable.
\end{ex}

\subsection{Other definitions} \label{sec:other}

There is a lot of other (equivalent) definitions of amenability of a countable
group, which illustrates importance and naturalness of this notion. We shall
briefly mention just some of them, referring the reader to \cite{Greenleaf69},
\cite{Pier84} and \cite{Paterson88} for more details. Moreover, the notion of
amenability has been extended to objects other than groups, in particular, to
group actions, equivalence relations, and, more generally, to groupoids (e.g.,
see \cite{Anantharaman-Renault00}).

The main application of the notion of amenability is its characterization as a
\emph{fixed point property}. Namely, \emph{a countable group $G$ is amenable
if and only if any continuous affine action of $G$ on a compact space has a
fixed point}. An example of such an action arises in the following way. Let
$X$ be a compact topological space endowed with a continuous action of $G$,
and let $\Pr(X)$ denote the space of probability measures on $X$ endowed with
the weak$^*$ topology. Then $\Pr(X)$ has a natural affine structure, and the
action of $G$ extends to a continuous affine action on $\Pr(X)$. Therefore,
\emph{any continuous action of an amenable group on a compact space has a
finite invariant measure}\footnotemark (in fact, this property can also be
shown to be equivalent to amenability). In the case of the group of integers
$\Z$ this result is known as the \emph{Krylov--Bogolyubov theorem}
\cite{Krylov-Bogolyubov37} (which is one of the starting points of the modern
theory of topological dynamical systems).

\footnotetext{The first proof of this fact by Bogolyubov \cite{Bogolyubov39}
published in 1939 (immediately after \cite{Krylov-Bogolyubov37}) in a rather
obscure journal in Ukrainian remained almost unknown, see
\cite{Anosov94,Ceccherini-Grigorchuk-delaHarpe99}.}

Yet another characterization of amenable groups can be given in terms of their
\emph{isoperimetric properties}. This condition is basically a specialization
of Reiter's condition to sequences of measures of special kind (although
historically it was introduced by F{\o}lner \cite{Folner55} some 10 years
before Reiter). Let $A_n$ be a sequence of finite subsets of $G$, and let
$\l_n$ be the associated uniform probability measures on $A_n$. Then Reiter's
condition for the sequence $\l_n$ is equivalent to the following condition on
the sets $A_n$:
$$
\frac{|gA_n \triangle A_n|}{|A_n|} \toto_{n\to\infty} 0 \qquad\forall\,g\in G
\;,
$$
where $\triangle$ denotes the symmetric difference of two sets, and $|A|$ is
the cardinality of a finite set $A$. A sequence of sets $A_n$ satisfying the
above condition is called a \emph{F{\o}lner sequence}, and the condition itself
is called \emph{F{\o}lner's condition}. Obviously, F{\o}lner's condition implies
Reiter's condition; however, the usual ``slicing'' isoperimetric techniques
also allow one to prove the converse, so that \emph{F{\o}lner's condition is
equivalent to amenability}.

For finitely generated groups F{\o}lner's condition takes especially simple form.
Indeed, in this case it is enough to verify it for the elements $g$ from a
finite generating set $K\subset G$ only. Let us assume that $K$ is symmetric,
and denote by $\G=\G(G,K)$ the (left) \emph{Cayley graph} of the group $G$
determined by $K$ (i.e., the vertex set is $G$, and the edges are of the form
$(g,kg)$ with $g\in G$ and $k\in K$)\footnote{Elsewhere in this paper we shall
always deal with the \emph{right} Cayley graphs. However, in order to keep the
notations consistent, here it is more convenient to consider the left Cayley
graphs.}. For a set $A\subset G$ denote by $\pt A = \pt_K A\subset A$ its
\emph{boundary} in the Cayley graph, i.e., the set of all points from $A$
which have a neighbor from the complement of $A$. Then a sequence of sets
$A_n\subset G$ is F{\o}lner if and only if
$$
\frac{|\pt A_n|}{|A_n|} \toto_{n\to\infty} 0 \;.
$$
Existence of a sequence of sets as above is an \emph{isoperimetric
characterization of amenability}.

\begin{ex}
For the group $\Z$ with the standard generating set $\{\pm 1\}$ the boundary
of the segment $A_n=\{0,1,2,\dots,n\}$ consists of two points $\{0,n\}$,
whereas $|A_n|\to\infty$, so that $\{A_n\}$ is a F{\o}lner sequence.
\end{ex}

\subsection{Elementary groups} \label{sec:elem}

The class of amenable groups is closed with respect to the ``elementary''
operations of taking subgroups, quotients, extensions and inductive limits.
Finite and abelian groups are amenable (cf. \exref{ex:cesaro}). The minimal
class of groups containing finite and abelian groups and closed with respect
to the above elementary operations is called \emph{elementary amenable} (EA)
\cite{Day57}.

In the above paper Day asked the question whether every amenable group is
elementary amenable. The first example of an amenable but not elementary
amenable group is the group of intermediate growth (see below) found by
Grigorchuk~\cite{Grigorchuk80,Grigorchuk85}. Later, a finitely presented
amenable extension of the Grigorchuk group was constructed
in~\cite{Grigorchuk98}.

However, there is yet another way to obtain ``obviously amenable'' groups. It
is related with the notion of \emph{growth}. Let $G$ be a finitely generated
group with a symmetric generating set $K$. Denote by $B_n$ the $n$-ball of the
Cayley graph metric on $G$ centered at the group identity, or, in other words,
the set of all elements of $G$ which can be presented as products of not more
than $n$ generators from $K$. The sequence $|B_n|$ is submultiplicative, so
that there exists a limit $\lim\log|B_n|/n$. The group $G$ is said to have
\emph{exponential} or \emph{subexponential growth} depending on whether this
limit is positive or zero (this property does not depend on the choice of a
generating set $K$).

The class of groups of subexponential growth contains all the groups of
\emph{polynomial growth} (the ones for which $|B_n|$ is bounded from above by
a polynomial function; by a theorem of Gromov \cite{Gromov81} these are
precisely finite extensions of nilpotent groups), but there are also examples
of groups of \emph{intermediate growth}, i.e., the ones whose growth is
subexponential without being polynomial. First examples of this kind were
constructed by Grigorchuk \cite{Grigorchuk85}, and these groups can often be
realized as self-similar groups (see
\cite{Bartholdi-Grigorchuk-Nekrashevych03}). For instance, the most famous of
the Grigorchuk groups has 4 generators acting on the rooted binary tree with
the matrix presentation
$$
 a \mapsto \begin{pmatrix} 0 & 1 \\ 1 & 0 \end{pmatrix} \;, \quad
 b \mapsto \begin{pmatrix} a & 0 \\ 0 & c \end{pmatrix} \;, \quad
 c \mapsto \begin{pmatrix} a & 0 \\ 0 & d \end{pmatrix} \;, \quad
 d \mapsto \begin{pmatrix} 1 & 0 \\ 0 & b \end{pmatrix} \;.
$$

If the growth of $G$ is subexponential, then the sequence $B_n$ necessarily
contains a F{\o}lner subsequence, so that \emph{the groups of subexponential
growth are amenable}. Therefore, one can change the definition of elementary
amenable groups by extending the set of ``building blocks'': the minimal class
of groups containing finite, abelian and subexponential groups and closed with
respect to the elementary operations is called \emph{subexponentially
elementary amenable} (SEA). Thus, a natural goal is to find amenable groups
which are not subexponentially elementary
(see~\cite{Grigorchuk98,Ceccherini-Grigorchuk-delaHarpe99}).

It is the Basilica group $\BB$ which provided the first example of this kind.
It was shown in~\cite{Grigorchuk-Zuk02a} that it does not belong to the class
SEA, whereas it was proved in~\cite{Bartholdi-Virag05} that the Basilica group
is amenable. We shall now explain the role of random walks in the proof of
amenability of $\BB$.

\begin{rem}
It is worth mentioning at this point that \emph{non-amenable groups do exist},
and there is actually quite a lot of them (for instance, numerous matrix
groups: non-elementary Fuchsian and Kleinian groups, lattices in semi-simple
Lie groups, etc.). The first example is of course the \emph{free group} $\F_d$
with $d\ge 2$ generators (since any discrete group can be obtained from free
groups by the elementary operations above, amenability of free groups would
have implied amenability of \emph{all} groups). Non-amenability of $\F_d$ may
be explained in many different ways by using various definitions. For
instance, it is not hard to build a \emph{paradoxical decomposition} of $\F_d$
(see \cite{Ceccherini-Grigorchuk-delaHarpe99} and the references therein),
which prevents it from having an invariant mean. Another way consists in
noticing that there are continuous actions of $\F_d$ on compact sets admitting
no finite invariant measures. Let, for instance $d=2$. Then an action of
$\F_2$ on a compact $K$ is determined by specifying two homeomorphisms of $K$
corresponding to the generators of $\F_2$. If $K$ is the circle, then the only
measure preserved by any irrational rotation is the Lebesgue measure. Take
such a rotation for the first homeomorphism from the definition of the action.
Then, if we choose the second homeomorphism in such a way that it does not
preserve the Lebesgue measure, then these two homeomorphisms (therefore, the
associated action of $\F_2$) do not have any invariant measure, so that $\F_2$
is not amenable.
\end{rem}

\section{Random walks} \label{sec:3}

\subsection{Convolution powers}

The main idea behind the use of random walks for establishing amenability at
first glance looks counterproductive. Let us replace arbitrary approximatively
invariant sequences of probability measures $\l_n$ on $G$ from Reiter's
condition with sequences of very special form, namely, with sequences of
convolution powers $\mu^{*n}$ of a single probability measure $\mu$ on $G$. In
the same way as with F{\o}lner's characterization of amenability (see
\secref{sec:other}), it turns out that this restricted class of
approximatively invariant sequences is still sufficient in order to
characterize amenability. More precisely, \emph{a group $G$ is amenable if and
only if there exists a probability measure $\mu$ on $G$ such that the sequence
of its convolution powers $\mu^{*n}$ is approximatively invariant.} In one
direction (the one we need for proving amenability) this is just a particular
case of Reiter's condition, whereas in the opposite direction it was
conjectured by Furstenberg \cite{Furstenberg73} and later independently proved
by Kaimanovich--Vershik \cite{Vershik-Kaimanovich79,Kaimanovich-Vershik83} and
by Rosenblatt \cite{Rosenblatt81}.

Now, working with sequences of convolution powers instead of arbitrary
sequences of probability measures on $G$ is actually easier, because of their
description as one-dimensional distributions of \emph{random walks} on $G$.
Moreover, one can use a very powerful quantitative criterion of whether a
given sequence of convolution powers is approximatively invariant. It is
provided by the \emph{entropy theory of random walks} which we shall briefly
describe below.

\medskip

Random walks are in a sense the most homogeneous Markov chains: they are
homogeneous both in space and in time. The latter property means that the
assignment $x\mapsto\pi_x$ of transition probabilities to the points from the
state space is equivariant with respect to a group action, and the simplest
instance of such an action is, of course, the action of a group on itself.

More formally, the (right) \emph{random walk} on a countable group $G$
determined by a probability measure $\mu$ is the Markov chain with the state
space $G$ and the transition probabilities
$$
\pi_g(g') = \mu(g^{-1}g')
$$
equivariant with respect to the left action of the group on itself. In other
words, from a point $g$ the random walk moves at the next moment of time to
the point $gh$, where the random \emph{increment} $h$ is chosen according to
the distribution $\mu$. We shall use for this description of transition
probabilities of the random walk $(G,\mu)$ the notation
$$
g \mapstoto_{h\sim\mu} gh \;.
$$
Thus, if the random walk starts at moment 0 from a point $g_0$, then its
position at time $n$ is
$$
g_n = g_0 h_1 h_2 \dots h_n \;,
$$
where $h_i$ is a \emph{Bernoulli sequence} of independent $\mu$-distributed
increments. Therefore, the distribution of the position $g_n$ of this random
walk at time $n$ is the translate $g_0\mu^{*n}$ of the \emph{$n$-fold
convolution power} $\mu^{*n}$ of the measure $\mu$ (the \emph{convolution} of
two probability measures $\mu_1,\mu_2$ on $G$ is defined as the image of the
product measure $\mu_1\otimes\mu_2$ on $G\times G$ under the map
$(g_1,g_2)\mapsto g_1g_2$).

\subsection{Trivial future}

Reiter's condition
$$
\|g\mu^{*n}-\mu^{*n}\|\toto_{n\to\infty} 0 \qquad \forall\,g\in G
$$
for the sequence of convolution powers $\mu^{*n}$ means that the time $n$
one-dimensional distributions of the random walk issued from the group
identity $e$ and from an arbitrary point $g\in G$ asymptotically coincide.
Probabilistically, the fact that the one-dimensional distributions of a Markov
chain asymptotically do not depend on their starting points means that the
``remote future'' (behaviour at infinity) of the chain does not depend on its
present, which, in view of the classical Kolmogorov argument means that the
``remote future'' must be trivial.

\begin{figure}[htp]
\centering \includegraphics[scale=.8]{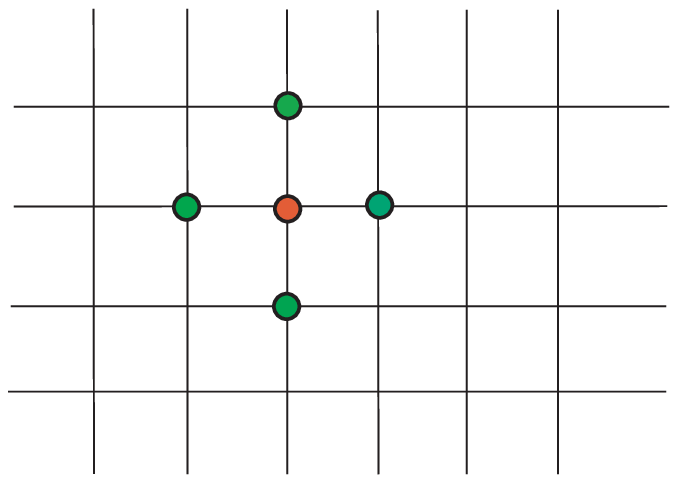} \hskip 3cm
\includegraphics[scale=0.46]{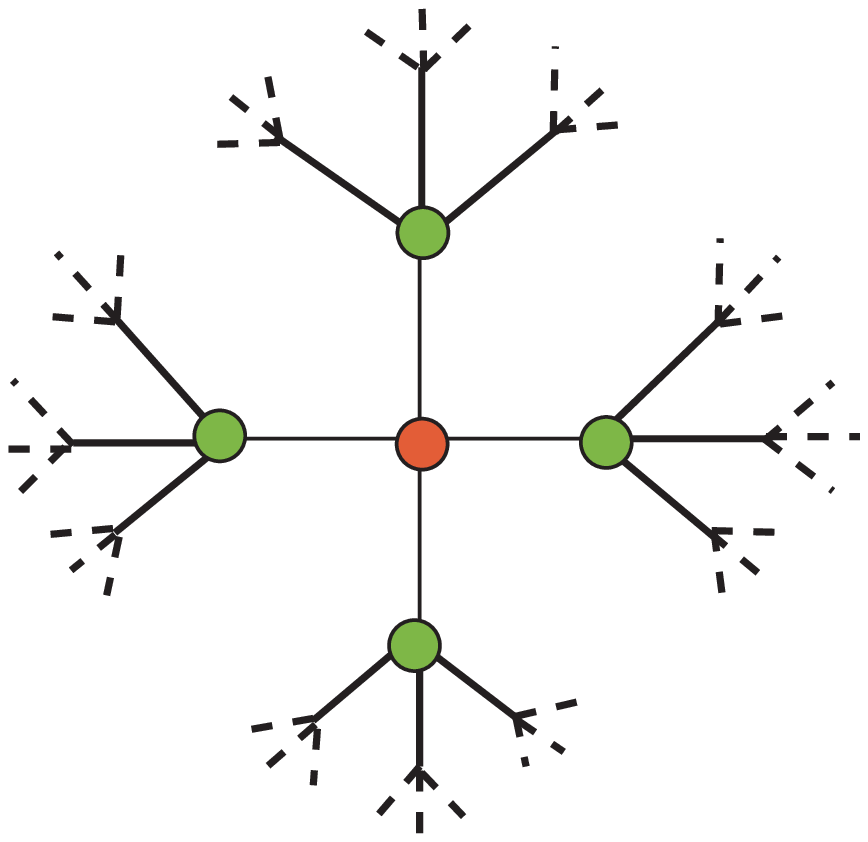}
\caption{} \label{fig:graphs}
\end{figure}

In order to explain the latter notion in a more rigorous way, let us look at
two examples: the simple random walks on the 2-dimensional integer lattice
($\equiv$ the Cayley graph of the group $\Z^2$) and on the homogeneous tree of
degree 4 ($\equiv$ the Cayley graph of the free group $\F_2$ with 2
generators). The simple random walk on a graph is the one whose transition
probabilities are equidistributed among neighbours; the simple random walk on
the Cayley graph of a group is precisely the random walk on this group
determined by the probability measure $\mu=\mu_K$ equidistributed on the
generating set $K$. Locally, each of these graphs is regular of degree 4 (each
point has 4 neighbours), but their global geometry is very different, see
\figref{fig:graphs}. In particular, the Cayley graph of $\F_2$ is endowed with
a natural boundary $\pt\F_2$ which can, for instance, be identified with the
space of all paths without backtracking issued from the identity of the group
(or any other reference point).

\begin{figure}[h]
\begin{center}
     \psfrag{d}[cl][cl]{${}$}
     \psfrag{z}[cl][cl]{${}$}
          \includegraphics[scale=.9]{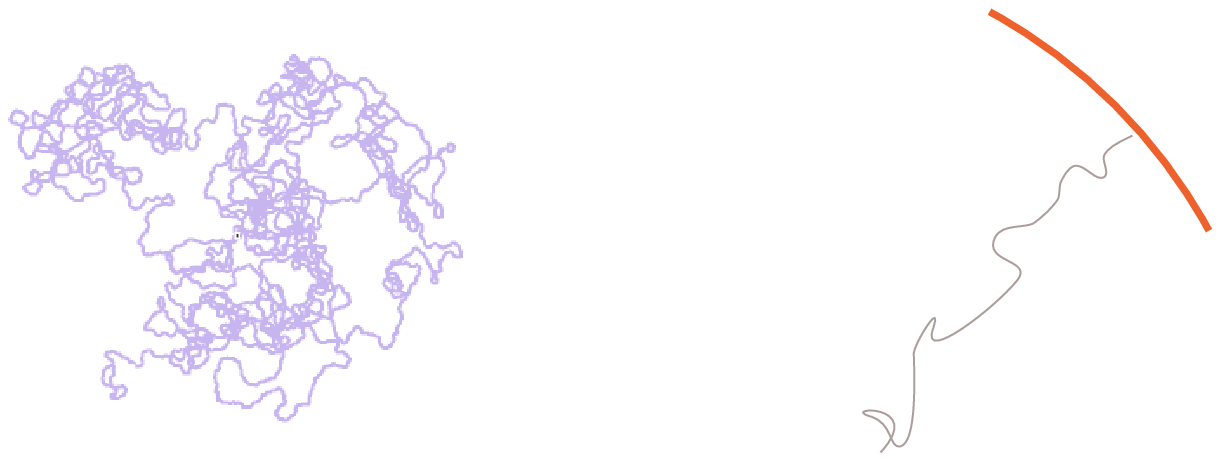}
          \end{center}
          \caption{}
          \label{fig:future}
\end{figure}

The global behaviour of sample paths of simple random walks on these graphs is
also very different. Sample paths of the simple random walk on $\F_2$ converge
a.s.\ to the boundary $\pt\F_2$ (of course, different sample paths may
converges to different limits). Thus, these limit points can be used to
distinguish sample paths by their behaviour at infinity. On the other hand,
although sample paths of the random walk on $\Z^2$ are quite complicated
(their scaling limit is the Brownian motion on the plane), they all look ``the
same'', see \figref{fig:future}.

Formally speaking, the ``remote future'' of a Markov chain is described by its
\emph{tail $\s$-algebra}
$$
\A^\infty= \bigcap_n \A_n^\infty \;,
$$
which is the limit of the decreasing sequence of $\s$-algebras $\A_n^\infty$
generated by the positions of the chain at times $\ge n$. Thus, the boundary
convergence of sample paths of the simple random walk on $\F_2$ at once
implies non-triviality of its tail $\s$-algebra. On the other hand, in spite
of absence of any visible behavior at infinity, proving triviality of the tail
$\s$-algebra for the simple random walk on $\Z^2$ requires additional work.

On a formal level a criterion of the triviality of the tail $\s$-algebra of an
arbitrary Markov chain is provided by the corresponding \emph{0--2 law}
\cite{Derriennic76,Kaimanovich92}. Its ``zero part'' for the random walk on a
group $G$ determined by a probability measure $\mu$ takes the following form
(see \cite{Kaimanovich-Vershik83}): the tail $\s$-algebra of the random walk
is trivial if and only if
$$
\| g\mu^{*n} - \mu^{*(n+1)} \| \toto_{n\to\infty} 0 \qquad
\forall\,g\in\supp\mu \;.
$$
Thus, \emph{if the support of the measure $\mu$ generates $G$ as a group, and
the tail $\s$-algebra of the random walk $(G,\mu)$ is trivial, then} the
sequence of Cesaro averages of the convolution powers satisfies Reiter's
condition, and therefore \emph{the group $G$ is amenable} (actually, if $\mu$
is \emph{aperiodic}, then the sequence of convolution powers $\mu^{*n}$ also
satisfies Reiter's condition).

\subsection{Asymptotic entropy}

Let $p=(p_i)$ be a discrete probability distribution. By Shannon the
\emph{amount of information} contained in any outcome $i$ is $-\log p_i$. The
average amount of information per outcome
$$
H(p) = - \sum p_i \log p_i
$$
is called the \emph{entropy} of the distribution $p$. The entropy has the
following two properties (by which it is essentially characterized, see
\cite{Yaglom-Yaglom83}):
\begin{itemize}
    \item \emph{monotonocity}: if $p'$ is a quotient of the distribution $p$
    (i.e., $p'$ is obtained by gluing together some of the states of $p$), then
    $$
    H(p')<H(p)\;,
    $$
    \item \emph{additivity}: for any two distributions $p_1,p_2$ the entropy
    of their product is
    $$
    H(p_1\otimes p_2) = H(p_1) + H(p_2)\;.
    $$
\end{itemize}

Let $\mu$ be a probability measure on a countable group $G$ with finite
entropy $H(\mu)$. The (asymptotic) \emph{entropy of the random walk} $(G,\mu)$
is defined as
$$
h(G,\mu) = \lim_{n\to\infty} \frac{H(\mu^{*n})}n \;,
$$
in other words, $h(G,\mu)$ is the asymptotic mean specific amount of
information about a single increment $h_i$ in the product $g_n=h_1h_2\dots
h_n$. Existence of the above limit follows from the fact that the sequence
$H(\mu^{*n})$ is subadditive. Indeed, by the definition of convolution, the
measure $\mu^{*(n+m)}$ is a quotient of the product measure
$\mu^{*n}\otimes\mu^{*m}$, whence
$$
H(\mu^{*(n+m)}) \le H(\mu^{*n}\otimes\mu^{*m}) = H(\mu^{*n}) + H(\mu^{*m})
$$
by the above monotonicity and additivity properties.

The asymptotic entropy was first introduced by Avez \cite{Avez72}. As it
turned out, it plays a crucial role in understanding the asymptotic properties
of the random walk $(G,\mu)$. Namely, as it was independently proved by
Kaimanovich--Vershik \cite{Vershik-Kaimanovich79,Kaimanovich-Vershik83} and
Derriennic \cite{Derriennic80}, \emph{the asymptotic entropy $h(G,\mu)$
vanishes if and only if the tail $\s$-algebra of the random walk $(G,\mu)$ is
trivial}.

This result means that the ``remote future'' of the random walk $(G,\mu)$ is
trivial if and only if the amount of information about the first increment
$h_1$ in the random product $g_n=h_1 h_2 \dots h_n$ asymptotically vanishes.
Indeed, on a more rigorous level, this amount of information is the \emph{mean
conditional entropy} (see \cite{Rohlin67}) $H\left(h_1|g_n\right)$ of $h_1$
with respect to $g_n$. It can be easily seen to coincide with the difference
$H\left(\mu^{*n}\right) - H\left(\mu^{*(n-1)}\right)$ between the entropies of
two consecutive convolution powers, whence the claim.

\subsection{Self-similarity and RWIDF}

A \emph{random walk with internal degrees of freedom} (RWIDF) on a group $G$
is a Markov chain whose state space is the product of $G$ by a certain space
$X$ (the \emph{space of degrees of freedom}), and its transition probabilities
are equivariant with respect to the action of $G$ on itself. Thus, the
transition probabilities are
$$
(g,x) \mapstoto^{\mu_{xy}(h)} (gh,y) \;,
$$
where $M=\{\mu_{xy}:x,y\in X\}$ is a $|X|\times|X|$ matrix of subprobability
measures on $G$ such that
$$
\sum_y \| \mu_{xy} \| = 1 \qquad \forall\,x\in X \;.
$$

The image of the RWIDF $(G,M)$ under the map $(g,x)\mapsto x$ is the quotient
Markov chain on $X$ with the transition matrix
$$
P=(p_{xy}) \;, \qquad p_{xy}=\|\mu_{xy}\| \;,
$$
which is the image of the matrix $M$ under the augmentation map.

\medskip

There is a natural link between random walks on self-similar groups and random
walks with internal degrees of freedom. Let $\mu$ be a probability measure on
a self-similar group $G$. Then the associated random walk is
$$
g\mapstoto^{\mu(h)} gh \;.
$$
By using the self-similar embedding $G\hookrightarrow\sym(X;G)$ it gives rise
to the random walk on the generalized permutation group $\sym(X;G)$ with the
transition probabilities
$$
M^g \mapstoto^{\mu(h)} M^g M^h \;.
$$
Since multiplication of the matrix $M^g\in\sym(X;G)$ by the increment $M^h$ is
done row by row, we obtain the following Markov chain on the space of these
rows:
$$
R \mapstoto^{\mu(h)} R M^h \;.
$$
Due to the definition of the group $\sym(X;G)$ the rows of the corresponding
matrices can be identified with points of the product space $G\times X$ (each
row has precisely one non-zero entry, so that it is completely described by
the value of this entry and by its position). Therefore, the latter Markov
chain can be interpreted as a Markov chain on $G\times X$ whose transition
probabilities are easily seen to be invariant with respect to the left action
of $G$ on $G\times X$, i.e., as a random walk on $G$ with internal degrees of
freedom parameterized by the alphabet $X$. This RWIDF is described by the
transition probabilities matrix
$$
M^\mu = \left( \mu_{xy} \right) = \sum_g \mu(g) M^g \;.
$$

Recall that stopping a Markov chain at the times when it visits a certain
recurrent subset of the state space produces a new Markov chain on this
recurrent subset (the \emph{trace} of the original Markov chain). For a random
walk with internal degrees of freedom on $G\times X$ determined by a matrix
$M$ take for such a subset the copy $G\times\{x\}$ of the group $G$ obtained
by fixing a value $x\in X$. It is recurrent provided the quotient chain on $X$
is irreducible. The transition probabilities of the induced chain on
$G\times\{x\}$ are obviously equivariant with respect to the left action of
the group $G$ on itself (because the original RWIDF also has this property).
Therefore, the induced chain on $G\times\{x\}$ is actually the usual random
walk determined by a certain probability measure $\mu^x$ on $G$.

The measures $\mu^x,\,x\in X$ can be expressed in terms of the matrix $M$ as
$$
\begin{aligned}
\mu^x &= \mu_{xx} + M_{x \ov x} \left(I + M_{\ov x \ov x} + M^2_{\ov x \ov x}
+ \dots\right) M_{\ov x x} \\ &= \mu_{xx} + M_{x \ov x} \left(I-M_{\ov x \ov
x}\right)^{-1} M_{\ov x x} \;,
\end{aligned}
$$
where $M_{x\ov x}$ (resp., $M_{\ov x x}$) denotes the row $(\mu_{xy})_{y\neq
x}$ (resp., the column $(\mu_{yx})_{y\neq x}$) of the matrix $M$ with the
removed element $\mu_{xx}$, and $M_{\ov x \ov x}$ is the $(d-1)\times (d-1)$
matrix (where $d=|X|$) obtained from $M$ by removing its $x$-th row and
column. The multiplication above is understood in the usual matrix
sense.\footnotemark

\footnotetext{As it is pointed out in \cite{Grigorchuk-Nekrashevych07}, this
formula corresponds to the classical operation of taking the \emph{Schur
complement} of a matrix.}

This is elementary probability. We look at the quotient chain on $X$ and
replace its transition probabilities $p_{xy}$ with the transition measures
$\mu_{xy}$ in the identity
$$
\begin{aligned}
1 &= p_{xx} + \sum_{n=0}^{\infty} \sum_{y_0,\dots,y_n\neq x} p_{x y_0} p_{y_0
y_1} \cdots p_{y_{n-1} y_n} p_{y_n x} \\ &= p_{xx}  + P_{x \ov x} \left(I +
P_{\ov x \ov x} + P^2_{\ov x \ov x} + \dots\right) P_{\ov x x} \\ &= p_{xx} +
P_{x \ov x} \left(I-P_{\ov x \ov x}\right)^{-1} P_{\ov x x} \;,
\end{aligned}
$$
which yields
$$
\begin{aligned}
\mu^x &= \mu_{xx} + \sum_{n=0}^{\infty} \sum_{y_0,\dots,y_n\neq x} \mu_{x y_0}
\mu_{y_0 y_1} \cdots \mu_{y_{n-1} y_n} \mu_{y_n x} \\ &= \mu_{xx}  + M_{x \ov
x} \left(I + M_{\ov x \ov x} + M^2_{\ov x \ov x} + \dots\right) M_{\ov x x} \\
&= \mu_{xx} + M_{x \ov x} \left(I-M_{\ov x \ov x}\right)^{-1} M_{\ov x x} \;.
\end{aligned}
$$
The first term in this formula corresponds to staying at the point $x$ (and
performing on $G$ the jump determined by the measure $\mu_{xx}$), whereas in
the second term the first factor corresponds to moving from $x$ to
$X\setminus\{x\}$, the second one to staying in $X\setminus\{x\}$ (each matrix
power $M^n_{\ov x \ov x}$ corresponding to staying in $X\setminus\{x\}$ for
precisely $n$ steps), and the third one to moving back from $X\setminus\{x\}$
to the point $x$. The matrix notation automatically takes care of what is
going on with the $G$-component of the RWIDF.

The measures $\mu^x$ also admit the following interpretation in terms of the
original random walk $(G,\mu)$ with the sample paths $(g_n)$: we look at it
only at the moments $n$ when $g_n(T_x)=T_x$, and $\mu^x$ is then the law of
the induced random walk on the group $\Aut(T_x)\cong \Aut(T)$.

\subsection{The M\"{u}nchhausen trick}

Let us now look at what happens with the asymptotic entropy in the course of
the transformations
$$
\mu \mapsto M=M^\mu \mapsto \mu^x \;.
$$
Since the information contained in a matrix does not exceed the sum of
informations about each row, passing to asymptotic entropies we obtain the
inequality
$$
h(G,\mu) \le d h(G,M) \;.
$$
Further, all points $x\in X$ are visited by RWIDF $(G,M)$ with the same
asymptotic frequency $1/d$ (because the uniform distribution is stationary for
the quotient chain on $X$), whence
$$
h(G,\mu^x)= d h(G,M) \ge h(G,\mu) \;.
$$

The idea of the \emph{M\"{u}nchhausen trick}\footnotemark \cite{Kaimanovich05}
consists in combining two observations. The first one is the above entropy
inequality. The second observation is that if the measure $\mu^x$ is a
non-trivial convex combination
$$
\mu^x = (1-\a)\d_e + \a\mu \;, \qquad 0<\a<1
$$
of the original measure $\mu$ and the $\d$-measure at the identity of the
group (in which case we call the measure $\mu$ \emph{self-similar}), then
$$
h(G,\mu^x) = \a h(G,\mu) \;.
$$
The result of these two observations is the inequality
$$
h(G,\mu) \le \a h(G,\mu) \;.
$$
Taken into account that $0<\a<1$, it is only possible if $h(G,\mu)=0$, which
proves amenability of the group $G$.

\footnotetext{It is named after venerable Baron M\"{u}nchhausen who
\textit{``\dots once rode on a cannon-ball and next told about it. Another
time he reported he had to get himself and the good horse he sat on, out of a
quagmire by pulling his own hair till he saved himself and his horse \dots''}}

We shall now show how the M\"{u}nchhausen trick works in two particular cases: for
the Basilica group $\BB$ and for the Mother groups $\mthr$.

\subsection{The Basilica group}

As we have seen in \secref{sec:bas}, the Basilica group is defined by the
matrix presentation
$$
\BB: a \mapsto \begin{pmatrix} b & 0 \\ 0 & 1 \end{pmatrix} \;, \qquad
 b \mapsto \begin{pmatrix} 0 & a \\ 1 & 0 \end{pmatrix} \;.
$$
Take on $\BB$ a symmetric probability measure $\mu$ supported by the
generators $a,b$ and their inverses\footnotemark:
$$
\mu = \a \left( a + a^{-1} \right) +\b \left( b + b^{-1} \right) \;,
$$
for $\a,\b>0$ such that $2(\a+\b)=1$. The matrix associated with the measure
$\mu$ is
$$
\begin{aligned}
M^\mu &= \a M^a + \a M^{a^{-1}} + \b M^b + \b M^{b^{-1}} \\
&= \a \begin{pmatrix} b & 0 \\ 0 & 1 \end{pmatrix} + \a \begin{pmatrix} b^{-1}
& 0 \\ 0 & 1 \end{pmatrix} + \b \begin{pmatrix} 0 & a \\ 1 & 0
\end{pmatrix} + \b \begin{pmatrix} 0 & 1 \\ a^{-1} & 0 \end{pmatrix} \\
&= \begin{pmatrix} \a(b+b^{-1}) & \b(1+a) \\ \b(1+a^{-1}) & 2\a
\end{pmatrix} \;.
\end{aligned}
$$
Therefore, the trace of the RWIDF $(\BB,M^\mu)$ on the copy of $\BB$
corresponding to the first letter 1 of the 2-letter alphabet $\{1,2\}$ is the
random walk governed by the measure
$$
\begin{aligned}
\wt\mu=\wt\mu^1 &= \mu_{11} + \mu_{12} \left( 1 - \mu_{22} \right)^{-1} \mu_{21} \\
&= \a(b+b^{-1}) + \frac{\b}2 (1+a)(1+a^{-1}) = \b + \frac{\b}2(a+a^{-1}) +
\a(b+b^{-1}) \;.
\end{aligned}
$$
If
$$
\frac\a\b = \frac\b{2\a} \iff 2\a^2 = \b^2 \;,
$$
then
$$
\wt\mu = \b + (1-\b) \mu \;,
$$
so that the measure $\mu$ is self-similar, and M\"{u}nchhausen's trick is
applicable.

\footnotetext{For simplicity we pass from now on to the group algebra
notations putting $g$ for a $\d$-measure $\d_g,\,g\neq e$ and just 1 for the
$\d$-measure $\d_e$ concentrated at the identity of the group.}

\subsection{The Mother group}

For proving amenability of the Mother group $\mthr$ (here and below we omit
the alpabet $X$) one can apply an approach somewhat different from the one
which was used above for the Basilica group. It is based on the fact that the
Mother group $\mthr$ is generated by two finite subgroups $A$ and $B$. We take
for the measure $\mu$ the convolution product of the uniform measures $m_A$
and $m_B$ on these subgroups:
$$
\mu = m_A * m_B \;.
$$
Then the matrix $M^\mu$ has a very special form
$$
M^\mu = M^{\mu_A} M^{\mu_B} = E_d \begin{pmatrix} \mu_B & 0 \\ 0 & \mu_A
E_{d-1}
\end{pmatrix} \;,
$$
where $d=|X|$, and $E_d$ denotes the order $d$ matrix with entries $1/d$, so
that $M^\mu$ has identical rows with entries
$$
M^\mu_{xy} =
  \begin{cases}
    \mu_B/d  &  \text{if $y=o$} \;, \\
    \mu_A/d  & \text{otherwise} \;.
  \end{cases}
$$
It means that transition probabilities of the associated RWIDF $(\mthr\times
X, M^\mu)$ do not depend on $x$, so that its projection to $\mthr$ is just the
random walk $(\mthr,\wt\mu)$ determined by the measure
$$
\wt\mu = \sum_y M^\mu_{xy} = \frac{d-1}{d} \mu_A + \frac{1}{d} \mu_B \;,
$$
whereas the projection of RWIDF $(\mthr\times X, M^\mu)$ to $X$ is the
sequence of independent $X$-valued random variables with uniform distribution
on $X$ (because all entries $M^\mu_{xy}$ have mass $1/d$). Thus,
$$
h(\mthr,\mu) \le d\, h(\mthr,\wt\mu) \;.
$$

The measure $\wt\mu$ is a convex combination of the idempotent measures
$\mu_A$ and $\mu_B$, so that its convolution powers are essentially convex
combinations of the convolution powers of $\mu$. The total number of $m_B$'s
in the $n$-fold convolution of $\wt\mu$ is $\sim n/d$, but some of them
disappear because $m_B m_B = m_B$, so that in fact $\wt\mu^{*n}$ is the
convolution of about $\frac{d-1}{d^2}n$ copies of $\mu=m_A*m_B$. Thus,
$$
h(\mthr,\wt\mu) = \frac{d-1}{d^2} h(\mthr,\mu) \;,
$$
whence $h(\mthr,\mu)=0$.

\bibliographystyle{amsalpha}
\bibliography{D:/Sorted/MyTex/mine}

\providecommand{\bysame}{\leavevmode\hbox to3em{\hrulefill}\thinspace}
\providecommand{\MR}{\relax\ifhmode\unskip\space\fi MR }
\providecommand{\MRhref}[2]{%
  \href{http://www.ams.org/mathscinet-getitem?mr=#1}{#2}
}
\providecommand{\href}[2]{#2}
\begin{thebibliography}{CSGdlH99}

\bibitem[ADR00]{Anantharaman-Renault00}
Claire Anantharaman-Delaroche and Jean Renault, \emph{Amenable groupoids},
  Monographies de L'Enseignement Math\'ematique [Monographs of L'Enseignement
  Math\'ematique], vol.~36, L'Enseignement Math\'ematique, Geneva, 2000, With a
  foreword by Georges Skandalis and Appendix B by E. Germain. \MR{2001m:22005}

\bibitem[Ano94]{Anosov94}
D.~V. Anosov, \emph{On {N}. {N}. {B}ogolyubov's contribution to the theory of
  dynamical systems}, Uspekhi Mat. Nauk \textbf{49} (1994), no.~5(299), 5--20.
  \MR{1311227 (96a:01029)}

\bibitem[Ave72]{Avez72}
Andr{\'e} Avez, \emph{Entropie des groupes de type fini}, C. R. Acad. Sci.
  Paris S\'er. A-B \textbf{275} (1972), A1363--A1366. \MR{0324741 (48 \#3090)}

\bibitem[BG00]{Bartholdi-Grigorchuk00}
L.~Bartholdi and R.~I. Grigorchuk, \emph{On the spectrum of {H}ecke type
  operators related to some fractal groups}, Tr. Mat. Inst. Steklova
  \textbf{231} (2000), no.~Din. Sist., Avtom. i Beskon. Gruppy, 5--45.
  \MR{1841750 (2002d:37017)}

\bibitem[BG02]{Bartholdi-Grigorchuk02}
Laurent Bartholdi and Rostislav Grigorchuk, \emph{On a group associated to
  $z^2-1$}, arXiv: math.GR/0203244, 2002.

\bibitem[BGN03]{Bartholdi-Grigorchuk-Nekrashevych03}
Laurent Bartholdi, Rostislav Grigorchuk, and Volodymyr Nekrashevych, \emph{From
  fractal groups to fractal sets}, Fractals in {G}raz 2001, Trends Math.,
  Birkh\"auser, Basel, 2003, pp.~25--118. \MR{2091700 (2005h:20056)}

\bibitem[Bie90]{Bielefeld90}
Ben Bielefeld (ed.), \emph{Conformal dynamics problem list}, Institute for
  Mathematical Sciences preprint series, no. 90-1, SUNY Stony Brook, 1990,
  arXiv: math.DS/9201271.

\bibitem[BKN08]{Bartholdi-Kaimanovich-Nekrashevych08}
Laurent Bartholdi, Vadim Kaimanovich, and Volodymyr Nekrashevych, \emph{On
  amenability of automata groups}, arXiv:0802.2837, 2008.

\bibitem[BN03]{Bondarenko-Nekrashevych03}
E.~Bondarenko and V.~Nekrashevych, \emph{Post-critically finite self-similar
  groups}, Algebra Discrete Math. (2003), no.~4, 21--32. \MR{2070400
  (2005d:20041)}

\bibitem[BN06]{Bartholdi-Nekrashevych06}
Laurent Bartholdi and Volodymyr Nekrashevych, \emph{Thurston equivalence of
  topological polynomials}, Acta Math. \textbf{197} (2006), no.~1, 1--51.
  \MR{2285317 (2008c:37072)}

\bibitem[Bog39]{Bogolyubov39}
N.~N. Bogolyubov, \emph{On some ergodic properties of continuous transformation
  groups}, Nauch. Zap. Kiev Univ. Phys.-Mat. Sb. \textbf{4} (1939), no.~5,
  45--52, in Ukranian, also \emph{Selected works in mathematics}, Fizmatlit,
  Moscow, 2006, pp. 213--222 (in Russian).

\bibitem[BS98]{Brunner-Sidki98}
A.~M. Brunner and Said Sidki, \emph{The generation of {${\rm GL}(n,\mathbb Z)$}
  by finite state automata}, Internat. J. Algebra Comput. \textbf{8} (1998),
  no.~1, 127--139. \MR{1492064 (99f:20055)}

\bibitem[BV05]{Bartholdi-Virag05}
Laurent Bartholdi and B{\'a}lint Vir{\'a}g, \emph{Amenability via random
  walks}, Duke Math. J. \textbf{130} (2005), no.~1, 39--56. \MR{2176547
  (2006h:43001)}

\bibitem[CSGdlH99]{Ceccherini-Grigorchuk-delaHarpe99}
Tullio Ceccherini-Silberstein, Rostislav~I. Grigorchuk, and Pierre de~la Harpe,
  \emph{Amenability and paradoxical decompositions for pseudogroups and
  discrete metric spaces}, Proc. Steklov Inst. Math. \textbf{224} (1999),
  no.~1, 57--97. \MR{1721355 (2001h:43001)}

\bibitem[CW89]{Connes-Woods89}
A.~Connes and E.~J. Woods, \emph{Hyperfinite von {N}eumann algebras and
  {P}oisson boundaries of time dependent random walks}, Pacific J. Math.
  \textbf{137} (1989), no.~2, 225--243. \MR{90h:46100}

\bibitem[Day49]{Day49}
Mahlon~M. Day, \emph{Means on semigroups and groups}, Bull. Amer. Math. Soc.
  \textbf{55} (1949), 1054--1055, abstract 55--11-507.

\bibitem[Day57]{Day57}
\bysame, \emph{Amenable semigroups}, Illinois J. Math. \textbf{1} (1957),
  509--544. \MR{0092128 (19,1067c)}

\bibitem[Day83]{Day83}
\bysame, \emph{Citation classic - amenable semigroups}, Current Contents Phys.
  Chem. Earth (1983), no.~26, 18--18.

\bibitem[Der76]{Derriennic76}
Yves Derriennic, \emph{Lois ``z\'ero ou deux'' pour les processus de {M}arkov.
  {A}pplications aux marches al\'eatoires}, Ann. Inst. H. Poincar\'e Sect. B
  (N.S.) \textbf{12} (1976), no.~2, 111--129. \MR{54 \#11508}

\bibitem[Der80]{Derriennic80}
\bysame, \emph{Quelques applications du th\'eor\`eme ergodique sous-additif},
  Conference on Random Walks (Kleebach, 1979) (French), Ast\'erisque, vol.~74,
  Soc. Math. France, Paris, 1980, pp.~183--201, 4. \MR{588163 (82e:60013)}

\bibitem[dlH00]{delaHarpe00}
Pierre de~la Harpe, \emph{Topics in geometric group theory}, Chicago Lectures
  in Mathematics, University of Chicago Press, Chicago, IL, 2000. \MR{1786869
  (2001i:20081)}

\bibitem[Fal03]{Falconer03}
Kenneth Falconer, \emph{Fractal geometry}, second ed., John Wiley \& Sons Inc.,
  Hoboken, NJ, 2003, Mathematical foundations and applications. \MR{2118797
  (2006b:28001)}

\bibitem[F{\o}l55]{Folner55}
Erling F{\o}lner, \emph{On groups with full {B}anach mean value}, Math. Scand.
  \textbf{3} (1955), 243--254. \MR{0079220 (18,51f)}

\bibitem[Fur73]{Furstenberg73}
Harry Furstenberg, \emph{Boundary theory and stochastic processes on
  homogeneous spaces}, Harmonic analysis on homogeneous spaces (Proc. Sympos.
  Pure Math., Vol. XXVI, Williams Coll., Williamstown, Mass., 1972), Amer.
  Math. Soc., Providence, R.I., 1973, pp.~193--229. \MR{50 \#4815}

\bibitem[GM05]{Glasner-Mozes05}
Yair Glasner and Shahar Mozes, \emph{Automata and square complexes}, Geom.
  Dedicata \textbf{111} (2005), 43--64. \MR{2155175 (2006g:20112)}

\bibitem[GN07]{Grigorchuk-Nekrashevych07}
Rostislav Grigorchuk and Volodymyr Nekrashevych, \emph{Self-similar groups,
  operator algebras and {S}chur complement}, J. Mod. Dyn. \textbf{1} (2007),
  no.~3, 323--370. \MR{2318495 (2008e:46072)}

\bibitem[Gre69]{Greenleaf69}
Frederick~P. Greenleaf, \emph{Invariant means on topological groups and their
  applications}, Van Nostrand Mathematical Studies, No. 16, Van Nostrand
  Reinhold Co., New York, 1969. \MR{40 \#4776}

\bibitem[Gri80]{Grigorchuk80}
R.~I. Grigorchuk, \emph{On {B}urnside's problem on periodic groups},
  Funktsional. Anal. Appl. \textbf{14} (1980), no.~1, 41--43. \MR{565099
  (81m:20045)}

\bibitem[Gri85]{Grigorchuk85}
\bysame, \emph{Degrees of growth of finitely generated groups and the theory of
  invariant means}, Math. SSSR Izv. \textbf{25} (1985), no.~2, 259--300.
  \MR{764305 (86h:20041)}

\bibitem[Gri98]{Grigorchuk98}
\bysame, \emph{An example of a finitely presented amenable group that does not
  belong to the class {EG}}, Sb. Math. \textbf{189} (1998), no.~1-2, 75--95.
  \MR{1616436 (99b:20055)}

\bibitem[Gro81]{Gromov81}
Mikhael Gromov, \emph{Groups of polynomial growth and expanding maps}, Inst.
  Hautes \'Etudes Sci. Publ. Math. (1981), no.~53, 53--73. \MR{623534
  (83b:53041)}

\bibitem[G{\.Z}02a]{Grigorchuk-Zuk02a}
Rostislav~I. Grigorchuk and Andrzej {\.Z}uk, \emph{On a torsion-free weakly
  branch group defined by a three state automaton}, Internat. J. Algebra
  Comput. \textbf{12} (2002), no.~1-2, 223--246, International Conference on
  Geometric and Combinatorial Methods in Group Theory and Semigroup Theory
  (Lincoln, NE, 2000). \MR{2003c:20048}

\bibitem[G{\.Z}02b]{Grigorchuk-Zuk02b}
\bysame, \emph{Spectral properties of a torsion-free weakly branch group
  defined by a three state automaton}, Computational and statistical group
  theory (Las Vegas, NV/Hoboken, NJ, 2001), Contemp. Math., vol. 298, Amer.
  Math. Soc., Providence, RI, 2002, pp.~57--82. \MR{2003h:60011}

\bibitem[Kai92]{Kaimanovich92}
Vadim~A. Kaimanovich, \emph{Measure-theoretic boundaries of {M}arkov chains,
  {$0$}-{$2$} laws and entropy}, Harmonic analysis and discrete potential
  theory (Frascati, 1991), Plenum, New York, 1992, pp.~145--180. \MR{94h:60099}

\bibitem[Kai95]{Kaimanovich95}
\bysame, \emph{The {P}oisson boundary of covering {M}arkov operators}, Israel
  J. Math. \textbf{89} (1995), no.~1-3, 77--134. \MR{96k:60194}

\bibitem[Kai03]{Kaimanovich03a}
\bysame, \emph{Random walks on {S}ierpi\'nski graphs: hyperbolicity and
  stochastic homogenization}, Fractals in Graz 2001, Trends Math.,
  Birkh\"auser, Basel, 2003, pp.~145--183. \MR{2091703 (2005h:28022)}

\bibitem[Kai05]{Kaimanovich05}
\bysame, \emph{``{M}\"unchhausen trick'' and amenability of self-similar
  groups}, Internat. J. Algebra Comput. \textbf{15} (2005), no.~5-6, 907--937.
  \MR{2197814}

\bibitem[KB37]{Krylov-Bogolyubov37}
Nicolas Kryloff and Nicolas Bogoliouboff, \emph{La th\'eorie g\'en\'erale de la
  mesure dans son application \`a l'\'etude des syst\`emes dynamiques de la
  m\'ecanique non lin\'eaire}, Ann. of Math. (2) \textbf{38} (1937), no.~1,
  65--113. \MR{1503326}

\bibitem[Kig01]{Kigami02}
Jun Kigami, \emph{Analysis on fractals}, Cambridge Tracts in Mathematics, vol.
  143, Cambridge University Press, Cambridge, 2001. \MR{1840042 (2002c:28015)}

\bibitem[KS83]{Kramli-Szasz83}
Andr{\'a}s Kr{\'a}mli and Domokos Sz{\'a}sz, \emph{Random walks with internal
  degrees of freedom. {I}. {L}ocal limit theorems}, Z. Wahrsch. Verw. Gebiete
  \textbf{63} (1983), no.~1, 85--95. \MR{85f:60098}

\bibitem[KV83]{Kaimanovich-Vershik83}
V.~A. Kaimanovich and A.~M. Vershik, \emph{Random walks on discrete groups:
  boundary and entropy}, Ann. Probab. \textbf{11} (1983), no.~3, 457--490.
  \MR{85d:60024}

\bibitem[Nek05]{Nekrashevych05}
Volodymyr Nekrashevych, \emph{Self-similar groups}, Mathematical Surveys and
  Monographs, vol. 117, American Mathematical Society, Providence, RI, 2005.
  \MR{2162164 (2006e:20047)}

\bibitem[Nek08]{Nekrashevych08}
\bysame, \emph{Free subgroups in groups acting on rooted trees},
  arXiv:0802.2554, 2008.

\bibitem[NT08]{Nekrashevych-Teplyaev08}
Volodymir Nekrashevych and Alexander Teplyaev, \emph{Groups and analysis on
  fractals}, Analysis on Graphs and its Applications (Proc. Sympos. Pure Math.,
  Vol. 77), Amer. Math. Soc., Providence, R.I., 2008, pp.~143--180.

\bibitem[Pat88]{Paterson88}
Alan L.~T. Paterson, \emph{Amenability}, Mathematical Surveys and Monographs,
  vol.~29, American Mathematical Society, Providence, RI, 1988. \MR{90e:43001}

\bibitem[Pie84]{Pier84}
Jean-Paul Pier, \emph{Amenable locally compact groups}, Pure and Applied
  Mathematics, John Wiley \& Sons Inc., New York, 1984, A Wiley-Interscience
  Publication. \MR{86a:43001}

\bibitem[Rei65]{Reiter65}
H.~Reiter, \emph{On some properties of locally compact groups}, Nederl. Akad.
  Wetensch. Proc. Ser. A 68=Indag. Math. \textbf{27} (1965), 697--701.
  \MR{0194908 (33 \#3114)}

\bibitem[Roh67]{Rohlin67}
V.~A. Rohlin, \emph{Lectures on the entropy theory of transformations with
  invariant measure}, Uspehi Mat. Nauk \textbf{22} (1967), no.~5 (137), 3--56.
  \MR{0217258 (36 \#349)}

\bibitem[Ros81]{Rosenblatt81}
Joseph Rosenblatt, \emph{Ergodic and mixing random walks on locally compact
  groups}, Math. Ann. \textbf{257} (1981), no.~1, 31--42. \MR{83f:43002}

\bibitem[RT09]{Rogers-Teplyaev09}
Luke~G. Rogers and Alexander Teplyaev, \emph{Laplacians on the {B}asilica
  {J}ulia set}, Commun. Pure Appl. Anal. (2009), to appear.

\bibitem[Sid00]{Sidki00}
Said Sidki, \emph{Automorphisms of one-rooted trees: growth, circuit structure,
  and acyclicity}, J. Math. Sci. (New York) \textbf{100} (2000), no.~1,
  1925--1943, Algebra, 12. \MR{1774362 (2002g:05100)}

\bibitem[VK79]{Vershik-Kaimanovich79}
A.M. Vershik and V.~A Kaimanovich, \emph{Random walks on groups: boundary,
  entropy, uniform distribution}, Dokl. Akad. Nauk SSSR \textbf{249} (1979),
  no.~1, 15--18. \MR{553972 (81f:60098)}

\bibitem[vN29]{vonNeumann29}
John von Neumann, \emph{Zur allgemeinen {T}heorie des {M}a{\ss}es}, Fund. Math.
  \textbf{13} (1929), 73--116 and 333, also \emph{Collected works}, vol.\ I,
  pages 599--643.

\bibitem[VV07]{Vorobets-Vorobets07}
Mariya Vorobets and Yaroslav Vorobets, \emph{On a free group of transformations
  defined by an automaton}, Geom. Dedicata \textbf{124} (2007), 237--249.
  \MR{2318547 (2008i:20030)}

\bibitem[YY83]{Yaglom-Yaglom83}
A.~M. Yaglom and I.~M. Yaglom, \emph{Probability and information}, Theory and
  Decision Library, vol.~35, D. Reidel Publishing Co., Dordrecht, 1983,
  Translated from the third Russian edition by V. K. Jain. \MR{736349
  (85d:94001)}

\end{thebibliography}

\enddocument

\bye